\documentclass[11pt]{article}

\usepackage{booktabs}
\usepackage[margin=15pt,font=small,labelfont={bf,sf},justification=justified]{caption}
\usepackage[superscript,biblabel]{cite}
\usepackage{color}
\usepackage{float}
\usepackage{graphicx}
\usepackage{lscape}
\usepackage[bbgreekl]{mathbbol}
\usepackage{mathtools}
\usepackage{multibib}
\usepackage{multirow}
\usepackage{tabularx}
\usepackage{titleref}
\usepackage{url}
\usepackage{lipsum}
\usepackage{bm}
\usepackage{tikz}
\usepackage{subfig}
\usepackage{caption}

\newif\ifbembo
\newif\ifcharter
\newif\iferewhon
\newif\iflibertine
\newif\iflibertinealt
\newif\ifpalantino
\newif\iftimesnewroman

\bembofalse
\chartertrue
\erewhonfalse
\libertinefalse
\libertinealtfalse
\palantinofalse
\timesnewromanfalse
\ifbembo
  \usepackage[p,osf]{ETbb} % osf in text, tabular lining figures in math
  \usepackage[scaled=.95,type1]{cabin} % sans serif in style of Gill Sans
  \usepackage[varqu,varl]{zi4}% inconsolata typewriter
  \usepackage[T1]{fontenc}
  \usepackage[libertine,vvarbb]{newtxmath}
  \usepackage[cal=boondoxo,bb=boondox]{mathalfa}
\fi

\ifcharter
  \usepackage[scaled=.98]{XCharter}
  \usepackage[scaled=1.04,varqu,varl]{inconsolata}
  \usepackage[type1]{cabin}
  \usepackage[uprightscript,xcharter,vvarbb,scaled=1.05]{newtxmath}
  \usepackage[cal=euler,scr=boondoxo,bb=boondox]{mathalpha}  % prefer bb=bboldxLight when available
\fi

\iferewhon
  \usepackage[p,osf,scaled=.98,space]{erewhon} % scaling by .98 not really necessary
  \usepackage[varqu,varl]{inconsolata} % typewriter
  \usepackage[type1,scaled=.95]{cabin} % sans serif like Gill Sans
  \usepackage[utopia,vvarbb]{newtxmath}
\fi

\iflibertine
  \usepackage[sb]{libertinus}
  \usepackage[T1]{fontenc}
  \usepackage{textcomp}
  \usepackage[varqu,varl]{zi4}% inconsolata for mono, not LibertineMono
  \usepackage{libertinust1math} % slanted integrals, by default
  \usepackage[scr=boondoxo,bb=boondox]{mathalpha} %Omit bb=boondox for default libertinus bb
\fi

\iflibertinealt
  \usepackage[sb]{libertinus}
  \usepackage[T1]{fontenc}
  \usepackage{textcomp}
  \usepackage{libertinust1math} % slanted integrals, by default
  \usepackage[scr=boondoxo,bb=boondox]{mathalpha} %Omit bb=boondox for default libertinus bb
\fi

\ifpalantino
  \usepackage[largesc]{newpxtext}
  \usepackage[vvarbb]{newpxmath}
\fi

\iftimesnewroman
  \usepackage{newtxtext} %T1 is default encoding
  \usepackage[scaled=0.95]{inconsolata} % typewriter
  \usepackage[vvarbb]{newtxmath} % vvarbb gives STIX Bbb
\fi

\usepackage[T1]{fontenc}
\usepackage[final,protrusion=true,expansion=true]{microtype}

\usepackage{titling}
\usepackage{authblk} % note: titling needs to be loaded before authblk
\pretitle{\begin{center}\bfseries\sffamily\LARGE}
\posttitle{\end{center}}
\usepackage{abstract}
\usepackage{sectsty} % note: the main alternative package, titlesec, does not work with titleref
\allsectionsfont{\sffamily}

\usepackage[left=1in,right=1in,top=1in,bottom=1in,includefoot,heightrounded]{geometry}
\usepackage[hidelinks]{hyperref}

\newcites{supp}{SUPPLEMENTAL REFERENCES}

\makeatletter
\patchcmd{\LS@rot}{90}{-90}{}{}
\patchcmd{\endlandscape}{90}{-90}{}{}
\makeatother

\usepackage{algorithm}
\usepackage{algorithmicx}
\usepackage{algpseudocode}

\usepackage{color}
\usepackage{xspace}

\definecolor{vargreen}{rgb}{0.0, 0.5, 0.0}

\newcommand{\WING}{WING\xspace}

\newcommand{\euleriandx}{{\ensuremath{\Delta x}}}

\newcommand{\soliddomO}{{\Omega^{\text{s}}_0}}

\newcommand{\Chib}{\boldsymbol{\chi}}
\newcommand{\Fb}{\boldsymbol{F}}
\newcommand{\Ub}{\boldsymbol{U}}
\newcommand{\ub}{\boldsymbol{u}}
\newcommand{\Xb}{\boldsymbol{X}}
\newcommand{\psib}{\boldsymbol{\psi}}
\newcommand{\dXb}{\, \mathrm{d}\Xb}
\newcommand{\PP}{\mathbb{P}}
\newcommand{\Sigmab}{\mathbb{S}}
\newcommand{\xb}{\boldsymbol{x}}
\newcommand{\half}{{\ensuremath{\frac{1}{2}}}}
\newcommand{\xface}{{i + \half, j, k}}
\newcommand{\yface}{{i, j + \half, k}}
\newcommand{\zface}{{i, j, k + \half}}

\title{\WING: A Simple Windowed Nonorthogonalized Initial Guess Procedure for Repeated Matrix Solves}
\author[1,*]{David Wells}
\author[4]{Matthew G. Knepley}
\author[1,2,3]{Boyce E. Griffith}
\affil[1]{Department of Mathematics, University of North Carolina at Chapel Hill, Chapel Hill, NC, USA}
\affil[2]{Department of Biomedical Engineering, University of North Carolina at Chapel Hill, Chapel Hill, NC, USA}
\affil[3]{McAllister Heart Institute, University of North Carolina at Chapel Hill, Chapel Hill, NC, USA}
\affil[4]{Department of Computer Science and Engineering, University at Buffalo, State University of New York, Buffalo, NY, USA}
\affil[*]{Corresponding author, \texttt{drwells@email.unc.edu}}

\begin{document}
\maketitle

\begin{abstract}
\noindent
Many numerical methods require solution of a sequence of linear systems with the same matrix and similar right-hand sides.
Krylov subspace methods are a common tool for solving such linear systems, and a carefully chosen initial guess for the solution can reduce the total number of iterations, and thereby the total computational cost, required for convergence to a specified numerical tolerance.
This paper introduces the \WING algorithm, a modification of Fischer's second algorithm, which lowers the cost of forming an acceptably close initial guess by skipping orthogonalization and solving the possibly singular normal equations with a pseudoinverse.
We demonstrate the efficacy of the new algorithm, particularly for solving linear systems with coarse relative tolerances, with numerical benchmarks based on fluid-structure interaction, mantle convection, and earthquake models.
\end{abstract}

\section{Introduction}
\label{sec:introduction}
Consider the linear system
\begin{equation}
  A x = b
  \label{eq:linear-system}
\end{equation}
with given solution and right-hand side pairs $\{(x_1, b_1), (x_2, b_2), \cdots, (x_N, b_N)\}$.
Under special circumstances, such as solving a matrix equation resulting from a discretization of a system of ordinary differential equations (ODEs), the right-hand side vectors are relatively close in value to each-other.
If we approximate a solution to Equation~\eqref{eq:linear-system} with right-hand side $b_{N + 1}$ and an iterative method (i.e., a method that forms a sequence converging to $A^{-1} b_{N + 1}$) then we can leverage this property to define a starting point for this sequence, called the initial guess, $\tilde{x}$, that approximates $A^{-1} b_{N + 1}$.

Many initial guess routines rely on the observation that if
\begin{equation}
  b \approx \tilde{b} = \sum \alpha_i b_i
  \label{eq:interpolate-b}
\end{equation}
for yet to be determined coefficients $\{\alpha_i\}$ then we should choose
\begin{equation}
  \tilde{x} = \sum \alpha_i x_i
  \label{eq:interpolate-x}
\end{equation}
because
\begin{equation}
  A \sum \alpha_i x_i = \sum \alpha_i b_i = \tilde{b} \approx{b}.
  \label{eq:application-of-guess}
\end{equation}
Equations~\eqref{eq:interpolate-b}-\eqref{eq:application-of-guess} form the core of Fischer's foundational work on choosing initial guesses~\cite{fischer1998projection}.
An alternative discussed by Fischer is to solve with a modified right-hand side $b - b_1$, motivated by the observation that if the new right-hand side $b$ equals a previous right-hand side then the system is trivial~\cite{fischer1998projection}.
Typically, to aid in the calculation of the $\alpha$ coefficients, $b$ will be projected onto an orthogonalized space instead of the set of provided right-hand sides.
This substitution guarantees a unique set of coefficients, whereas the original set of right-hand sides $\{b_i\}$ may be linearly dependent.
In practice, this approach requires balancing the computational cost of reorthogonalization of a large set of vectors against the computational benefit of the improved initial guess.
In addition, initial guess algorithms based on projection typically reset the number of stored solution-right-hand side pairs to zero once they reach some preset capacity to control the overall computational cost and memory usage of the linear solver.
This limitation, resulting from the need for orthogonalization, was partially overcome by Christensen~\cite{christensen2017efficient} with a rolling QR factorization update algorithm, i.e., a method which can simultaneously add and remove vectors.
The primary drawback of the rolling QR method is that it requires relatively expensive updates to the stored vectors, which can increase the total computational cost.
Christensen modeled the cost of updating the projection space as approximately $t_{\text{p}} = 6 L N$, in which $6$ comes from work required to perform a Givens rotation, $N$ is the number of vectors, and $L$ is the length of each vector.
Analysis in Chapter 4 of Christensen~\cite{christensen2017efficient} concludes that the accuracy improvement in the initial guess, measured by the number of correct digits, scales proportionally to $\log(N)$ and that the most cost-effective search spaces sizes are small --- about six vectors.
The rolling QR method implemented by Christensen~\cite{christensen2017efficient} was extended by Austin et al.~\cite{austin2021initial}, which also compared a variety of other initial guess algorithms such as interpolation and least-squares extrapolation.
Section 6.3 of Austin et al.~\cite{austin2021initial} concludes that about eight vectors ``seems to be the dimension of diminishing returns'' for the projection methods.
L{\"o}hner~\cite{lohner2005projective} showed that, for certain problems with a small number of stored vectors, orthogonalization is not necessary.

Many alternative initial guess routines utilize properties of the underlying numerical discretization to achieve good results.
Ali and Sadkane~\cite{ali2012improved} use the linear structure of multi-step ODE solvers to form linear spaces that contain approximations to the solution.
Section 4 of Pitton et al.~\cite{pitton2019accelerating} proposed a similar method based on projecting onto an orthogonalized subspace.
Grinberg~\cite{grinberg2011extrapolation} proposed two methods based on extrapolation for accelerating linear solvers used in a Navier-Stokes solver.
Their first method is based on proper orthogonal decomposition (POD), a common technique for model reduction.
In that work, the initial guess is computed by a low-rank approximation of the system.
Grinberg's second proposed method is based on a local extrapolation technique that requires knowledge of the underlying physical discretization (in this case, high-order modal finite elements).
%% TODO'': `using the previous solution' is based on what is in Grinberg's
%% thesis, not the paper, which I cannot access right now
Both techniques have excellent performance and lower solver time by about 50\% relative to using the previous solution as the initial guess. %%% TODO: do we have a citation for 50\%? What does this mean in this context?

% WING: WIndowed Nonorthogonalized Guess

In contrast to these methods, this study proposes a new modification of Fischer's second algorithm~\cite{fischer1998projection}.
We introduce the \emph{WIndowed Nonorthogonalized Guess}, which we call the \WING algorithm.
This method avoids the potentially expensive orthogonalization step by simply skipping it and instead
solving the resulting (possibly singular) normal equations via least-squares.
Further, since this is a purely algebraic approach, the implementation of this method requires no problem-specific information and can be easily added to any existing iterative solver implementation.
One implementation is presently available in recent releases of PETSc~\cite{petsc-web-page} and can easily be enabled with any PETSc \texttt{KSP} solver.
Hence, this method combines the good accuracy properties of projection methods as discussed by Austin et al.~\cite{austin2021initial} while avoiding reorthogonalization, which is typically the most computationally expensive part of an initial guess procedure.

The primary contribution of this study is the introduction and analysis of the \WING algorithm.
We show that its cost is, among typical projection schemes, essentially optimal.
We demonstrate the utility of our algorithm with numerical experiments examining the relative efficacy of \WING compared to other initial guess algorithms.

% description of WING
\section{The \WING Algorithm}
\label{sec:algorithm}
The \WING algorithm defines the coefficients $\{\alpha_i\}$ in Equation~\eqref{eq:interpolate-x} as the least-squares projection coefficients of the new right-hand side onto the space spanned by previous right-hand sides.
Given $N$ previous solution and right-hand side pairs $\{(x_1, b_1), (x_2, b_2), \cdots, (x_N, b_N)\}$ and a new right-hand side $b$, the right-hand side correlation matrix $C$ and right-hand side vector $c$ are
\begin{equation}
  C_{i,j} = b_i \cdot b_j \text{ and } c_i = b_i \cdot b.
  \label{eq:define-least-squares-system}
\end{equation}
Because the right-hand side vectors may be linearly dependent, in practice the coefficients $\alpha = C^+ c$ are computed via the pseudoinverse $C^+$ of $C$, which is discussed in Section~\ref{sec:ill-conditioning}.
If all right-hand side vectors are linearly independent, then $C$ is invertible, and
\begin{equation}
  \tilde{b} = \sum_i \alpha_i b_i
  \label{eq:interpolate-b-2}
\end{equation}
is the orthogonal projection of $b$ onto the space spanned by the previous right-hand sides.
Hence, combining Equations~\eqref{eq:interpolate-x}, \eqref{eq:define-least-squares-system}, and \eqref{eq:interpolate-b-2} we obtain the definition of the \WING initial guess.

% 1. windowing, matrix-vector products
Like most initial guess algorithms that depend on previous solver states, the \WING algorithm consists of an \emph{update} routine and a \emph{formation} routine.
If the correlation matrix $C$ has $N$ rows and columns, then the update routine removes the oldest row and column.
Next, independent of the size of $C$, the update routine appends a new row and column computed from the most recent right-hand side.
Similarly, if there are $N$ stored vector pairs, then the oldest pair is removed and the new solution-right-hand side pair is added.
This is a windowed algorithm in the sense that, after some startup period, the spanning sets for the solution vector and right-hand side always contain the $N$ most recent such vectors.
In contrast, the two algorithms proposed by Fischer~\cite{fischer1998projection} and Pitton~\cite{pitton2019accelerating} are not windowed.
If those have $N$ vectors, then they remove all stored vectors, clear the correlation matrix, and start over with just the new solution-right-hand side pair.
In addition to discarding stored vectors, which will decrease the quality of the initial guess, orthogonalization with respect to the inner product induced by $A$ requires potentially expensive matrix-vector products.
For $|Sb|$ stored vectors, the projection algorithm of Pitton et al.~\cite{pitton2019accelerating} requires $|Sb|$ matrix-vector multiplications.
The \WING algorithm permits windowing because, unlike these algorithms, it does not build bases; instead, it keeps all provided vectors in the spanning sets independent of whether or not they are in the span of previously stored vectors.
Because no vectors are orthogonalized, no matrix-vector products are used, which further lowers the computational cost of \WING.
Consequently, unlike many projection methods that use an inner product induced by $A$, \WING does not require $A$ to be symmetric positive definite.

% 2. simplicity and cost
In addition to permitting windowing, a major advantage of \WING is its simplicity.
Unlike algorithms that recycle parts of the Krylov subspace~\cite{pitton2019accelerating} or exploit the structure of the linear system~\cite{ali2012improved}, \WING is a `black-box' solver that only requires solution and right-hand side vectors.
In particular, \WING does not require any knowledge of the linear system itself, just the vectors.
This makes it easy to implement.
This simplicity makes \WING computationally cheaper than standard projection methods.
One may reuse the dot products used to compute $c$ in Algorithm~\ref{alg:form} to update $C$ in Algorithm~\ref{alg:update}, thus lowering the total cost of the computation to approximately $L (2 N + 1) + N^3$ floating point operations corresponding to length $L$ vectors, $N + 1$ dot products (as the update routine must compute $b \cdot b$ in addition to the $N$ dot products necessary to compute $c$), $N$ scaled vector additions, and the least-squares solution of an $N \times N$ system.
As $N$ tends to be less than ten and thus $N \ll L$, we can ignore the cost of the pseudoinverse (and other $O(N^2)$ or $O(N^3)$ operations) when estimating the computational cost of \WING.
Indeed, the only $O(L)$ computational costs are the $N + 1$ dot products and $N$ scaled vector additions.
As typical projection schemes require computing an $N \times N$ correlation matrix and forming the guess with $N$ scaled vector multiplications, we say that \WING is essentially optimal because it does not incur any additional $O(L)$ costs.

% 3. use b, not A x
Another choice that improves \WING's performance is using the provided solution and right-hand side vectors from the converged linear solver.
Because the solution is computed only to a provided tolerance, in general $A x \neq b$ for the provided $x$ and $b$ pair (though, typically, $||A x - b|| \leq \varepsilon ||b||$ for some relative tolerance $\varepsilon$).
To avoid inconsistencies, many initial guess routines, such as the PETSc~\cite{petsc-web-page} implementation of Fischer's first and second methods, use $A x$ rather than $b$.
Because performing matrix-vector products is typically the most computationally costly step in forming the initial guess, we may be able to improve an algorithm's performance by skipping as many of them as possible.

% 4. summarize
The relative efficacy of \WING at strict and loose solver tolerances is investigated in Section~\ref{subsec:ibfe}.
Together, these choices make \WING's performance competitive with extrapolation methods that compute equivalent coefficients $c$ with a polynomial whose independent variable is time, despite \WING being a projection algorithm.
For example, an extrapolation method also using $N$ vectors will require $N$ scaled vector additions for a total cost of $L N$ --- i.e., half the computational cost of \WING.
Austin et al.~\cite{austin2021initial} found that even though projection required fewer solver iterations than extrapolation, the QR update routine was sufficiently expensive that extrapolation methods had overall better performance.

\begin{algorithm}
  \caption{\WING formation routine.}
  \label{alg:form}
  \begin{algorithmic}[1]
    \Require{$Sx$, $Sb$ are the spanning sets of previous solution and right-hand side vectors}
    \Require{$C$ is the right-hand side correlation matrix}
    \Require{$b$ is the current right-hand side}
    \Require{$\epsilon$ is the pseudoinverse tolerance}
    \Function{Formation}{$Sx$, $Sb$, $C$, $b$, $\epsilon$}
    \For{$i \gets 1 \textrm{ to } |Sb|$}
    \State $c(i) \gets Sb[i] \cdot b$
    \EndFor
    \State $\alpha \gets \mathrm{pseudoinverse}(C, \epsilon) c$
    \State $x \gets 0$
    \For{$i \gets 1 \textrm{ to } |Sx|$}
    \State $x \gets x + \alpha[i] Sx[i] $
    \EndFor
    \State \Return{$x$}
    \EndFunction
  \end{algorithmic}
\end{algorithm}

\begin{algorithm}
  \caption{\WING update routine.}
  \label{alg:update}
  \begin{algorithmic}[1]
    \Require{$Sx$, $Sb$ are the spanning sets of solution and right-hand side vectors}
    \Require{$x$, $b$ are the newest solution vector and right-hand side vector}
    \Require{$N$ is the maximum number of vector pairs}
    \Function{Update}{$Sx$, $Sb$, $x$, $b$, $N$}
    \If{$|Sx| = N$}
    \Comment{If we have the maximum number of vector pairs, remove the oldest one}
    \State $Sx$ $\gets$ $Sx$ - $\{\mathrm{oldest}(Sx)\}$ + $\{x\}$
    \State $Sb$ $\gets$ $Sb$ - $\{\mathrm{oldest}(Sb)\}$ + $\{b\}$
    \Else
    \State $Sx$ $\gets$ $Sx$ + $\{x\}$
    \State $Sb$ $\gets$ $Sb$ + $\{b\}$
    \EndIf
    \State $C \gets 0$
    \For{$i \gets 1 \textrm{ to } |Sb|$}
    \For{$j \gets 1 \textrm{ to } |Sb|$}
    \State $C(i, j) \gets Sb[i] \cdot Sb[j]$
    \EndFor
    \EndFor
    \State \Return{$\{Sx, Sb, C\}$}
    \EndFunction
  \end{algorithmic}
\end{algorithm}

\section{Overcoming ill-conditioning with the pseudoinverse}
\label{sec:ill-conditioning}
For \WING, the stored set of right-hand side vectors may not form a basis whereas the Fischer-2 algorithm \cite{fischer1998projection} only retains linearly independent right-hand side vectors.
Hence, a significant disadvantage of \WING is the increased difficulty in solving the resulting normal equations, which (as the set of right-hand side vectors is not orthonormalized) may require solving a singular or nearly singular matrix.
Our numerical experiments in Section~\ref{sec:numerical-results} indicate that the rank of the correlation matrix tends to be no more than $3$ or $4$, independent of the number of stored vectors.
Hence, in practice, we expect this linear system to be singular.
% Another viewpoint discussed in Section 2 of Austin et al.~\cite{austin2021initial} is that adding a new right-hand side vector to $Sb$ whose (after being projected onto $Sb$) residual is relatively small will only marginally improve the search space and may only increase the overall time to solution.
Regardless of its rank, we solve the system of normal equations in Algorithm~\ref{alg:form} with a Moore-Penrose pseudoinverse
\begin{align}
  C^+(\epsilon) &= V \Sigma^+ V^T,                                            \\
  \Sigma^+_{ij} &=
  \begin{cases}
    1 / \sigma_i, &\text{if } i = j \text{ and } \sigma_i > \epsilon \sigma_0,\\
    0 &\text{otherwise},
  \end{cases}
  \label{eq:pseudoinverse}
\end{align}
in which $\sigma_i$ is the $i$th largest singular value of $C$ and $V$ is the orthogonal matrix of right singular vectors (which, because $C$ is symmetric, is the transpose of the matrix of left singular vectors).
To see why this system has a solution, consider a set of vectors $Cb$ containing two duplicated vectors $b_i = b_j$.
Consequently $c_i = c_j$, so rows $i$ and $j$ of $c$ and also $C$ are equal.
Hence, applying Gaussian elimination to $C$ and $c$, we obtain a system in which the row-reduced version of $c$ will contain zeros in all rows not corresponding to pivot rows of $C$.
The same argument can be generalized to any linear combination, i.e., for any $r \leq N$ number of linearly independent vectors in $Cb$ the normal equations must have at least one solution.
Solving these normal equations with a pseudoinverse guarantees a unique solution by selecting the solution vector $c$ with the lowest Euclidean norm~\cite{demmel1997applied}.

\section{Numerical Results}
\label{sec:numerical-results}
\subsection{Initial Guess Algorithms}
In this section we compare four different initial guess algorithms, examining both their relative efficacy for evolution problems and more general linear systems in which we expect the right-hand side to change slowly.
% typesetting: its a list so don't capitalize 'the'
The relevant initial guess algorithms are the two projection methods described by Fischer~\cite{fischer1998projection}, a POD method from Zampini~\cite{zampini2010non}, and the \WING algorithm introduced herein.
Since Fischer's second algorithm requires a symmetric positive definite (SPD) linear operator, we do not use it with problems which do not solve SPD linear systems.
The POD-based method is fundamentally different from the other three because it uses information about the linear system and solution data to construct a reduced-order model (ROM).
Rather than relying on orthogonalization, the POD method uses the singular value decomposition to extract pertinent information from the sets of right-hand side and solution vectors.
The resulting ROM is then used to compute a predicted solution to the given linear system of equations.
As this method is also windowed it has nearly identical solver iteration counts to \WING, especially with less than four stored vectors.

All experiments were performed on a workstation.
The IBFE experiment in Section \ref{subsec:ibfe} and the mantle convection experiment in Section \ref{subsec:mantle-convection} were parallelized via MPI with four processes and the earthquake model experiment in Section \ref{subsec:pylith} used one MPI process.

\subsection{Immersed Boundary Finite Element / Difference Projection (IFED) Operators}
\label{subsec:ibfe}
This subsection investigates a fluid-structure interaction (FSI) problem consisting of a neutrally buoyant soft elastic ball interacting with an incompressible fluid.
We compute all results twice, using both a coarse ($10^{-8}$) and strict ($10^{-14}$) relative tolerance for both relevant linear solvers.
This example uses the immersed boundary finite element (IBFE) module of IBAMR~\cite{IBAMR} and is based on the lid-driven cavity example from Section 5.2 of Griffith and Luo~\cite{Griffith2017}.
This section only presents algorithmic details of this method relevant to the present study.
A complete description of the discretization is available in Griffith and Luo~\cite{Griffith2017}.

\subsubsection{Benchmark Description}
\begin{figure}
\centering
% scale is -0.25 to 0.0: double-check units!
\begin{tabular}{c   c   c}
\includegraphics[width=0.27\textwidth]{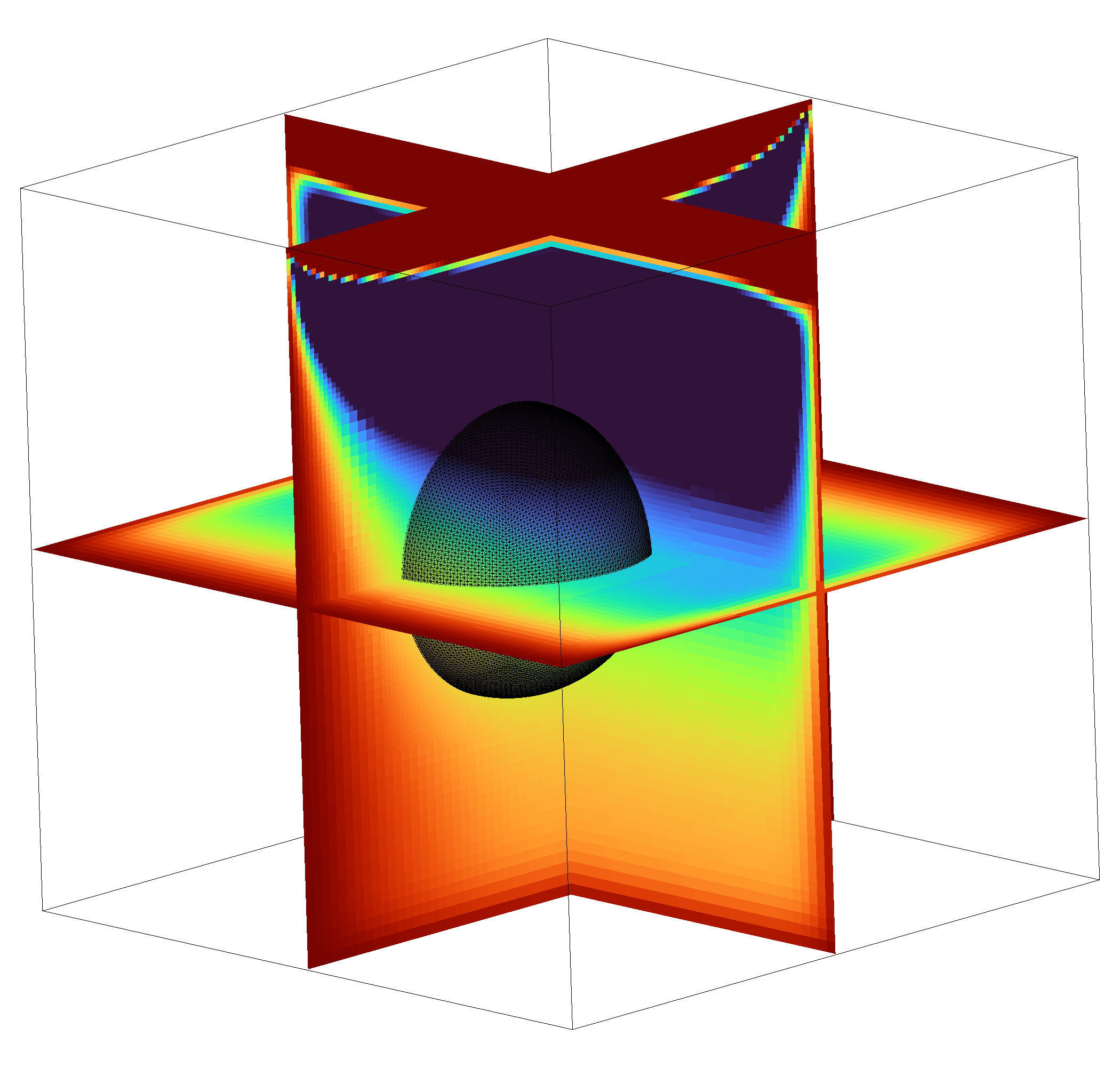} \hspace{0.05\textwidth} &
\includegraphics[width=0.27\textwidth]{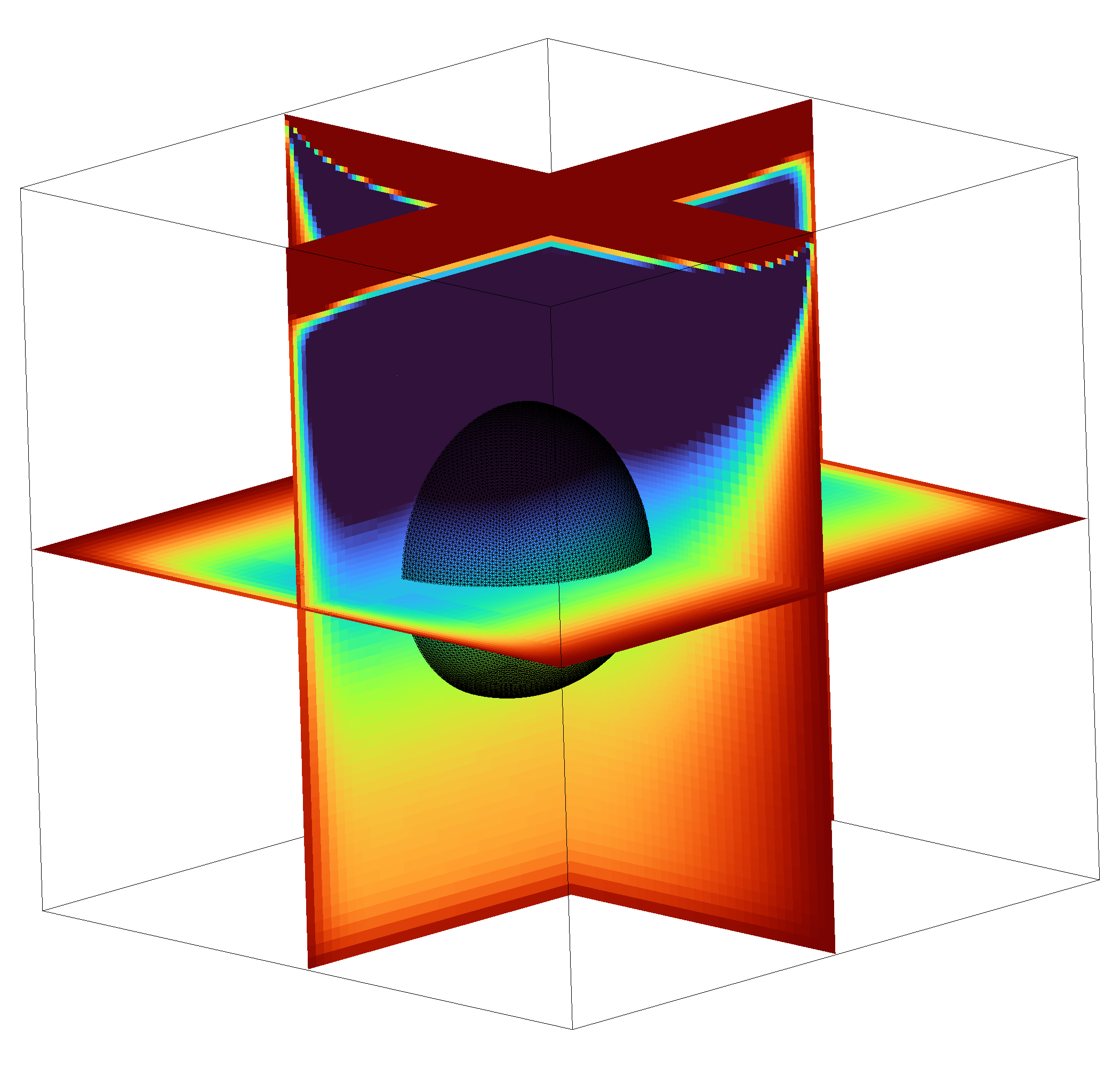} \hspace{0.05\textwidth} &
\includegraphics[width=0.27\textwidth]{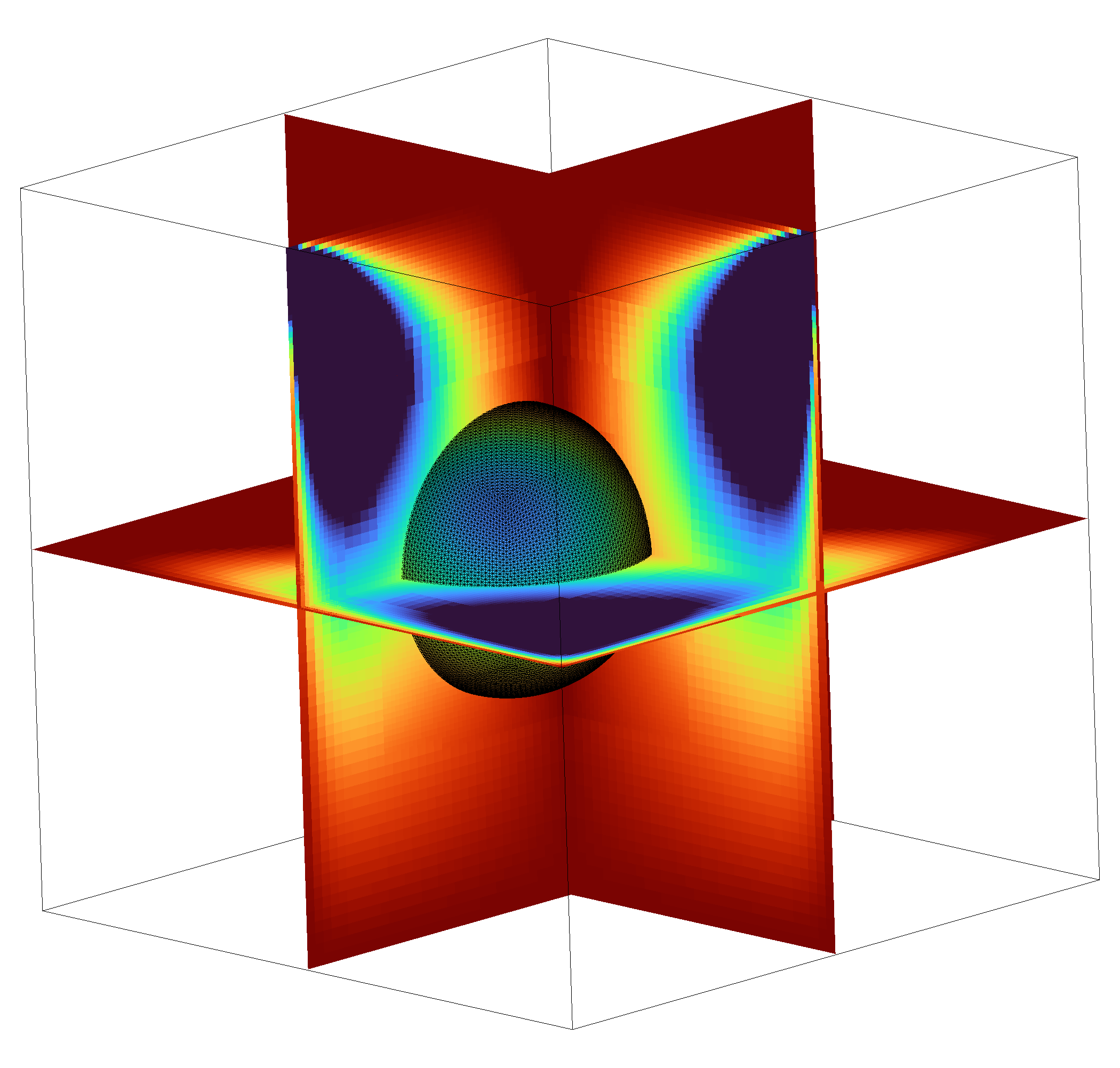}
\end{tabular}

% This spacing is eye-balled
\vspace{0.05\textwidth}
\includegraphics[width=0.4\textwidth]{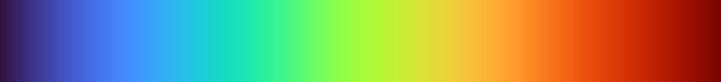}\\
$-0.25 \text{ cm/s}$\hspace{0.31\linewidth}$0 \text{ cm/s}$

\caption{$x$, $y$, and $z$ components of the Cartesian grid velocity for the IBFE-4 example after $500$ time steps.}
\end{figure}

In the IFED method, a structure's velocity, which is represented by a finite element (FE) field, is computed by interpolating and projecting a velocity defined on a Cartesian grid by
\begin{align}
    \Ub^{\text{IB},1}(\Xb,t) &= \sum_{i,j,k} u^1_{\xface}
    \, \delta_h(\xb_{\xface} - \Chib_h(\Xb,t)) \euleriandx^3
                                                                              \\
    \Ub^{\text{IB},2}(\Xb,t) &= \sum_{i,j,k} u^2_{\yface}
    \, \delta_h(\xb_{\yface} - \Chib_h(\Xb,t)) \euleriandx^3
                                                                              \\
    \label{eq:semidiscrete-interp-3}
    \Ub^{\text{IB},3}(\Xb,t) &= \sum_{i,j,k} u^3_{\zface}
    \, \delta_h(\xb_{\zface} - \Chib_h(\Xb,t)) \euleriandx^3                    \\
    \int_{\soliddomO} \Ub(\Xb,t) \cdot \psib(\Xb) \dXb &=
    \int_{\soliddomO} \Ub^{\text{IB}}(\Xb,t) \cdot \psib(\Xb) \dXb
\end{align}
in which $\Ub^{\text{IB}}$ is the velocity interpolated at selected points, $\Xb$ is a point on the structure in its reference configuration, $\xb$ is a point on the Cartesian grid, $\soliddomO$ is the reference configuration of the structure, $\Chib_h(\Xb,t)$ is the mapping from the reference to current configuration of the structure, $u^1, u^2$, and $u^3$ are the velocity components defined on the Cartesian grid, $\euleriandx$ is the Cartesian grid mesh spacing, $\delta_h(\xb)$ is a regularized delta function, $\psib(\Xb)$ is a test function, and $\Ub(\Xb)$ is the resulting FE velocity field.
The velocity values are defined on a staggered grid and evolved with the marker and cell (MAC) algorithm~\cite{harlow1965numerical}.
All finite element spaces use the standard 27-node hexahedral element (HEX27).

In this example the force is defined on a structural mesh discretized with the FE method.
In particular, the force is computed as projection of the weak divergence of the first Piola-Kirchoff stress $\PP$, which depends on the deformation gradient $\nabla_{\Xb} \Chib_h(\Xb, t)$ of the structure, onto the finite element space via
\begin{equation}
    \int_{\soliddomO} \Fb_h(\Xb,t) \cdot \psib(\Xb) \dXb =
    -\int_{\soliddomO} \PP(\Xb,t) : \nabla_{\Xb} \psib(\Xb) \dXb
    \label{eq:weak-force}
\end{equation}
in which $\soliddomO$ is the reference configuration of the structure, $\Xb \in \soliddomO$ is a point in the reference configuration, $\psib(\Xb)$ is a test function, and $\Fb_h(\Xb, t)$ is the resulting FE force field.

We use a timestepping algorithm which performs one force projection and two velocity projections at each time step, which requires solving three $L^2$ projections with slowly changing right-hand sides.
In contrast to the velocity, the evolution of the force is only implicitly dependent on the solution of a PDE.
Like the pressure, the force is not a state variable and is nonlinearly dependent on the deformation of the structure.
Hence, like what was observed for the pressure in Grinberg~\cite{grinberg2011extrapolation}, using a specialized initial guess routine does not lower the iteration counts for the force projection as much as the velocity projection.

The original Fischer algorithms perform very well with a large number of stored vectors.
However, since these methods do not use a windowing procedure to save old data, they provide less benefit with small numbers of stored vectors since the initial guess space is frequently cleared.
Hence, both the \WING and POD algorithms, which incrementally update the sets of stored vectors, show better performance if using relatively few vectors as compared to the Fischer algorithms for the same linear system.

Both the Fischer and POD algorithms require either orthogonalization or order reduction with the SVD.
Hence, these algorithms have higher computational costs compared to \WING, which uses the provided vectors to compute initial guesses directly without any postprocessing.
Because linearly dependent right-hand side vectors are not excluded from guess formation, the linear spaces stored by \WING tend to have lower rank than the spaces stored by the Fischer algorithms.
We conclude this is why the \WING iteration counts are either similar or worse than the Fischer iteration counts for higher numbers of stored vectors for both numerical tolerances.

\subsubsection{Solver performance with a relative tolerance of $10^{-8}$}
\label{subsubsec:ibfe-loose-tol}
\begin{figure}
  \centering
  \begin{tabular}{c c}
  \includegraphics[width=0.45\linewidth]{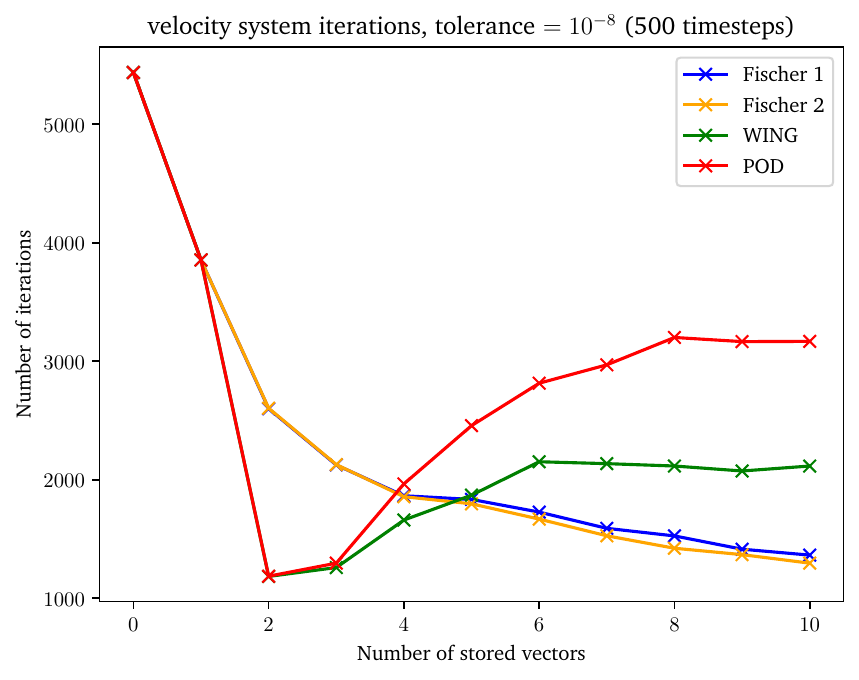} &
  \includegraphics[width=0.45\linewidth]{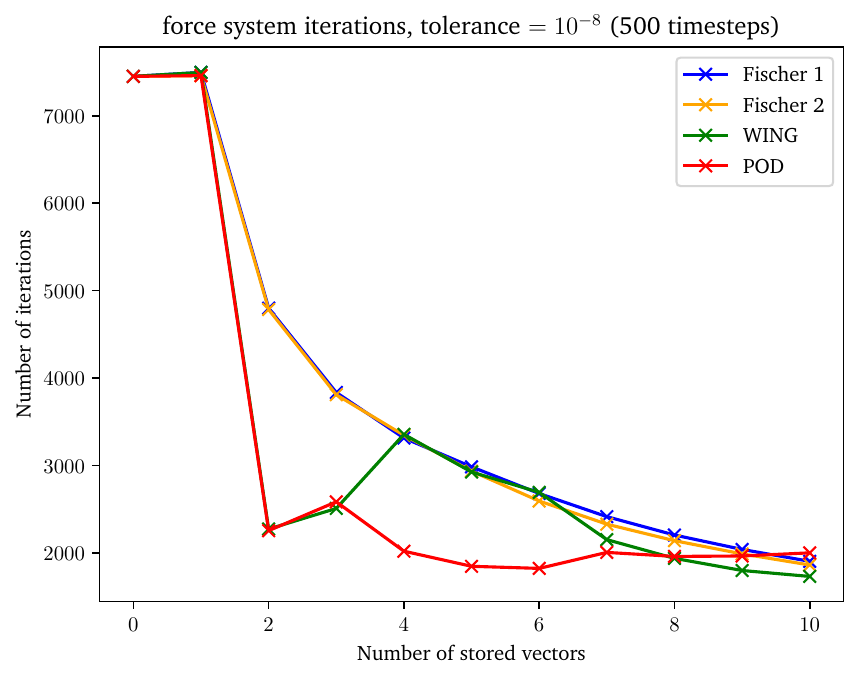}\\
  \end{tabular}

  \caption{Iteration counts for the IBFE-4 experiment for both finite element systems with a coarser relative tolerance.}

  \label{fig:ibfe-iterations-tol}
\end{figure}

\begin{figure}
  \centering
  \begin{tabular}{c c}
  \includegraphics[width=0.45\linewidth]{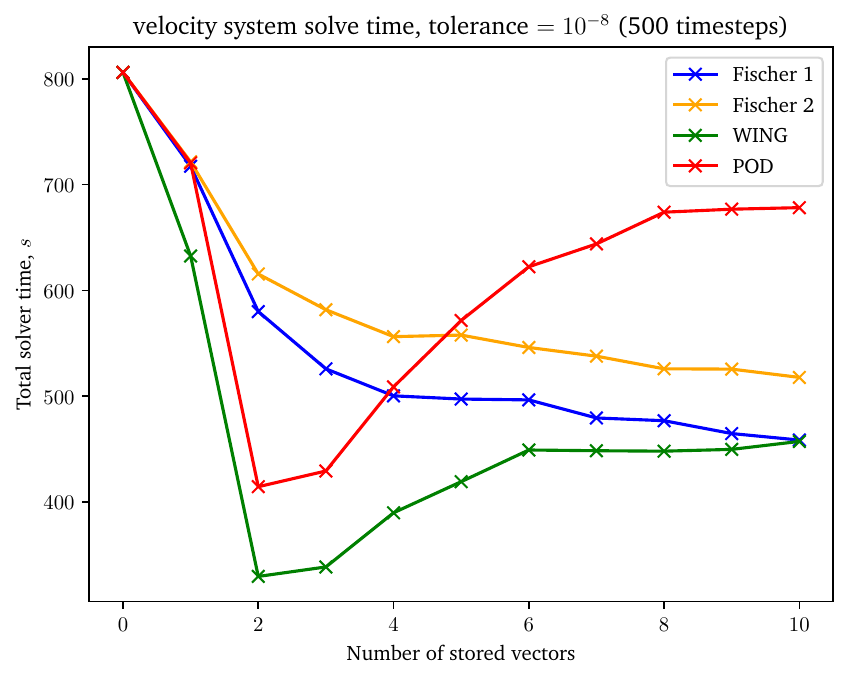} &
  \includegraphics[width=0.45\linewidth]{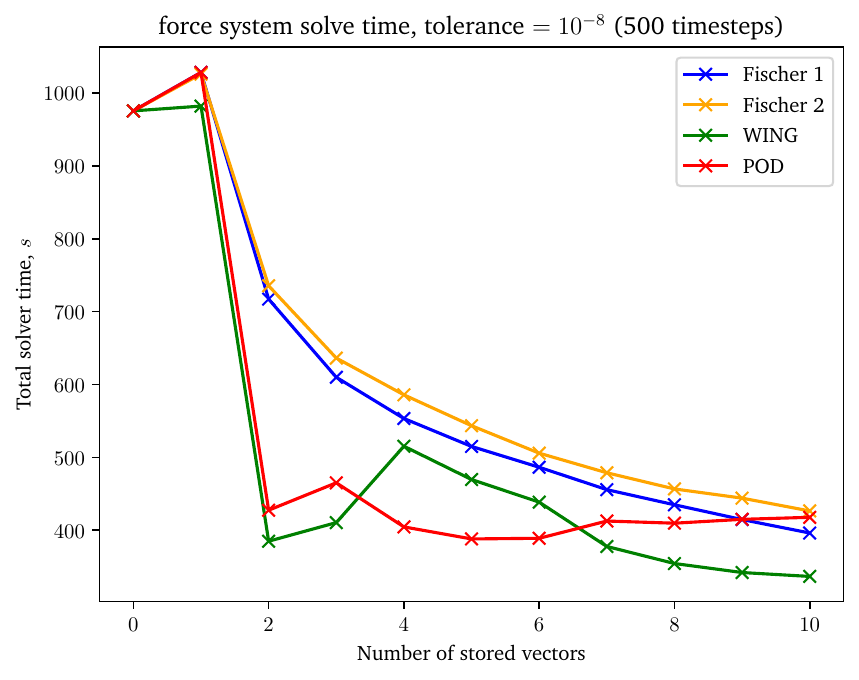}\\
  \end{tabular}

  \caption{Solver times for the IBFE-4 experiment for both finite element systems with a coarser relative tolerance.}

  \label{fig:ibfe-time-tol}
\end{figure}

\begin{figure}
  \centering
  \begin{tabular}{c c}
  \includegraphics[width=0.45\linewidth]{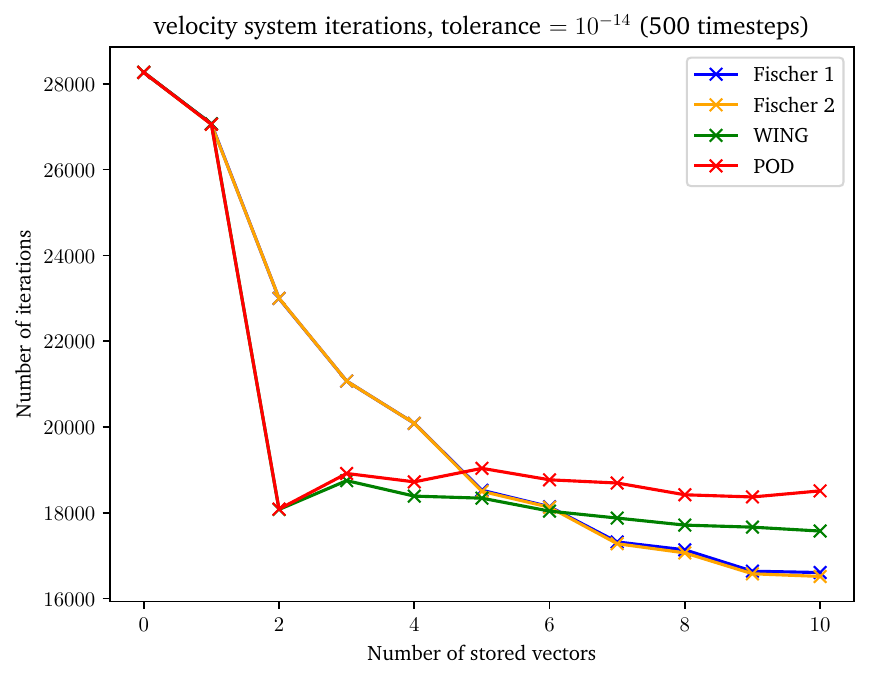} &
  \includegraphics[width=0.45\linewidth]{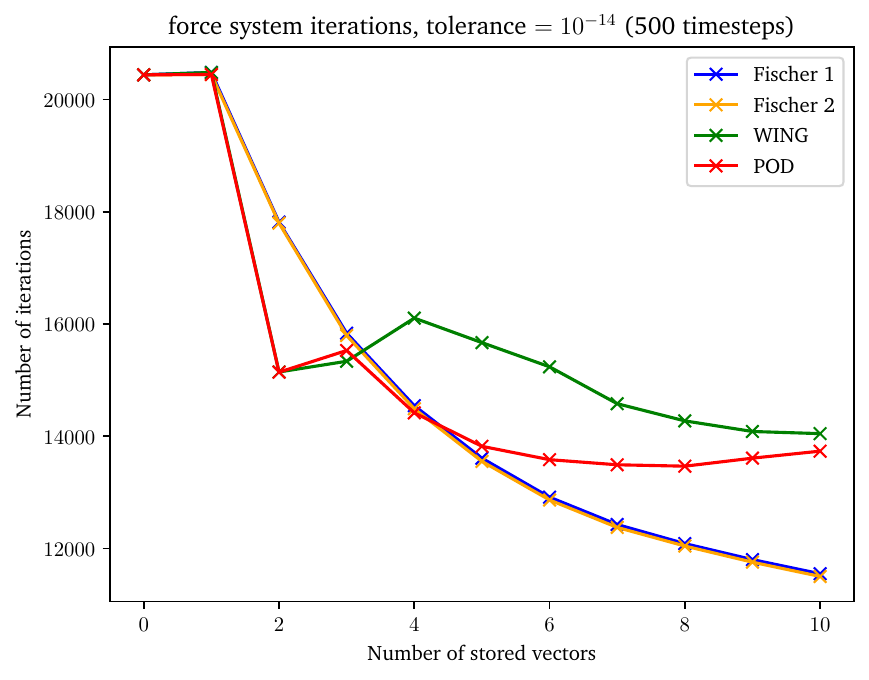}\\
  \end{tabular}

  \caption{Iteration counts for the IBFE-4 experiment for both finite element systems with a stricter relative tolerance.}

  \label{fig:ibfe-iterations}
\end{figure}

\begin{figure}
  \centering
  \begin{tabular}{c c}
  \includegraphics[width=0.45\linewidth]{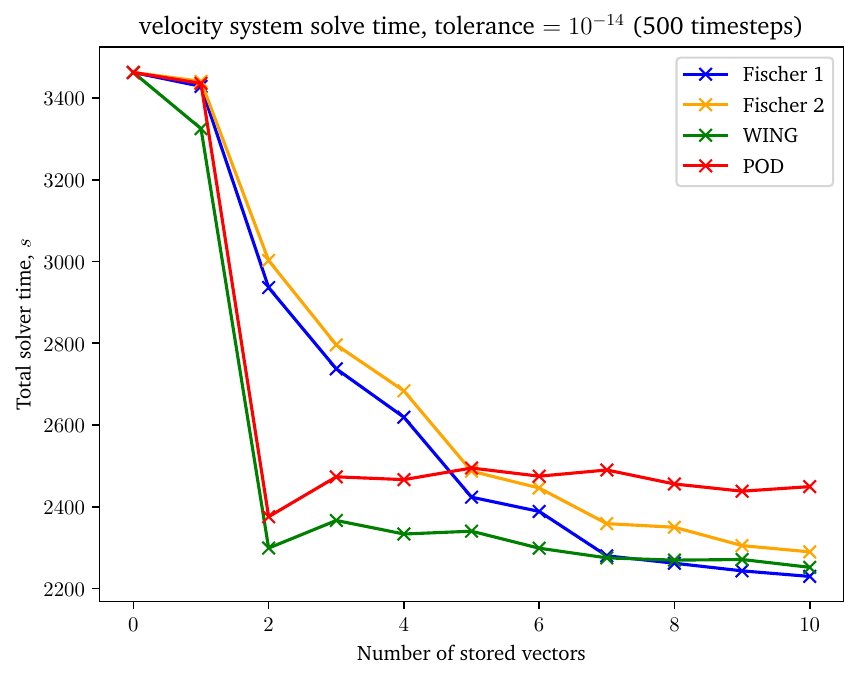} &
  \includegraphics[width=0.45\linewidth]{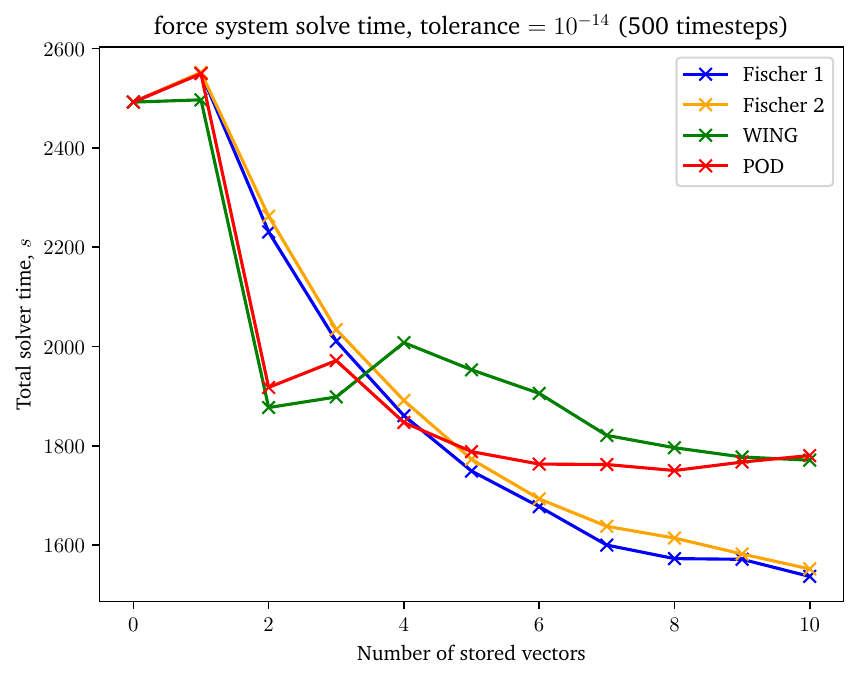}\\
  \end{tabular}

  \caption{Solver times for the IBFE-4 experiment for both finite element systems with a stricter relative tolerance.}

  \label{fig:ibfe-time}
\end{figure}

Figures~\ref{fig:ibfe-iterations} and \ref{fig:ibfe-time} demonstrate the runtime performance of the four initial guess routines for a variety of stored space sizes.
Both pairs of plots show the total number of iterations (summed over all time steps) as well as the total solver time for both velocity and force projection.
At this tolerance, the timestepping algorithm requires about one to three solver iterations for velocity projection (as two are performed per time step) and about four to $8$ solver iterations for force projection per time step.
The two original Fischer guess algorithms, because of orthogonalization, perform steadily better as the number of stored vectors increases.
In contrast, the POD and \WING algorithms have worse performance for velocity projection as more vectors are added.
This is likely a consequence of the lack of orthogonalization in both algorithms, i.e., the resulting sets of stored vectors may be linearly dependent, which may degrade the accuracy in the pseudoinverse step in Equation~\eqref{eq:pseudoinverse}.

Despite the higher iteration counts, the \WING algorithm still requires less total solver time than the original Fischer algorithms because it avoids doing any additional matrix-vector products.
Consequently, even though the stored set of vectors are linearly dependent and provide a worse initial guess, the \WING algorithm still requires less solution time due to its lower per-iteration computational requirements.
The \WING and POD algorithms also outperform the two Fischer algorithms when using just two or three stored vectors because they avoid the cost of frequent reorthogonalization.
Indeed, the sharpest improvements across all cases come from using no initial guess to using just two stored vectors for both \WING and POD.

\subsubsection{Solver performance with a relative tolerance of $10^{-14}$}
\label{subsubsec:ibfe-strict-tol}
In this subsection we examine the case in which the linear solves use a much stricter relative tolerance.
In this case the solver needs a much larger number of iterations to converge.
Here, the timestepping algorithm requires about $18$ to $20$ solver iterations for velocity projection (as two are performed per time step) and about $28$ to $32$ solver iterations for force projection per time step.
As a result, the choice of initial guess routine, while still important, has less of an impact on the total number of iterations than the $10^{-8}$ case.
For example, with the coarser tolerance, a good initial guess routine lowered the number of velocity iterations from $5438$ to $1364$.
In contrast, using the strict tolerance lowered the number of iterations from $28266$ to $17572$: i.e., from a factor of four improvement to less than a factor of two.

Given the stricter tolerance, the more computationally expensive initial guess routines performed best since their up-front costs were amortized over the relatively higher number of solver iterations.
Hence, the most efficient algorithms with respect to the number of total solver iterations are the two Fischer algorithms, each with a large number of stored vectors.
This is particularly true for the force projection, in which, like the coarser tolerance case, the methods with larger subspaces tend to require fewer overall solver iterations.
Hence, since the Fischer algorithms produce the highest quality spaces by orthogonalizing the input data, it is not surprising that they achieve the best performance.
Because the velocity projection does not benefit as much from larger guess spaces, \WING performs adequately well on solver time since it has lower computational requirements.

\subsection{PyLith}
\label{subsec:pylith}
\begin{figure}
  \centering
  \begin{tabular}{c c}
  \includegraphics[width=0.45\linewidth]{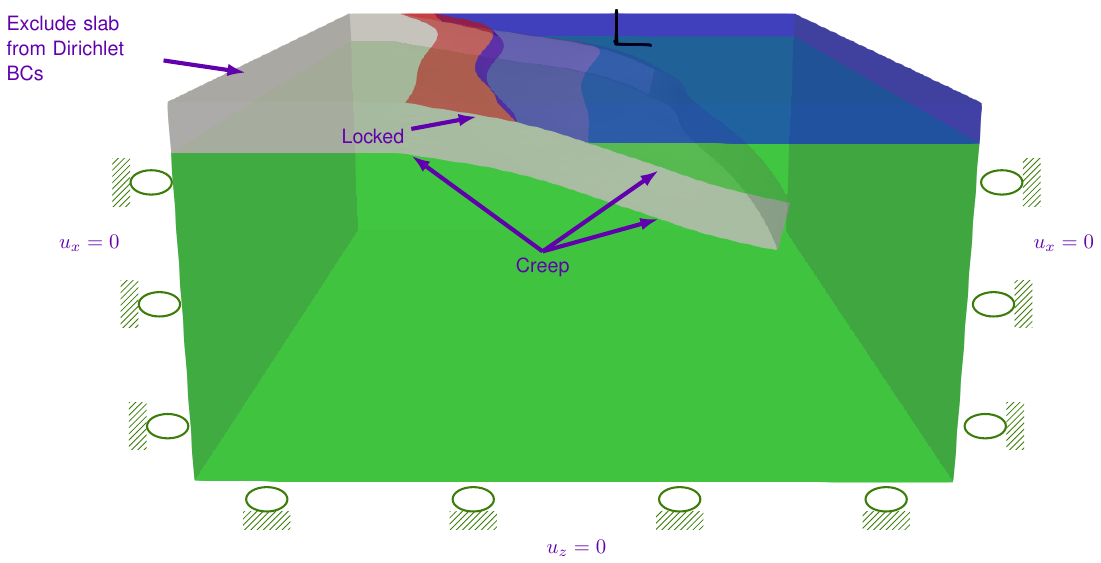}
  \includegraphics[width=0.45\linewidth]{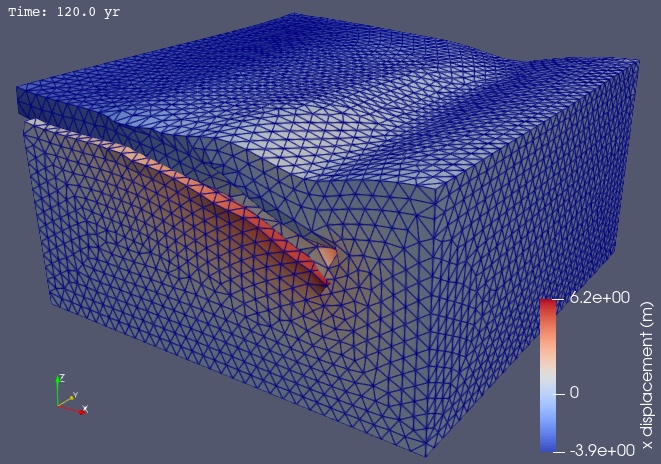}
  \end{tabular}

  \caption{Boundary conditions and representative solution for the PyLith subduction problem.}

  \label{fig:subduction-3d-step03-diagram}
\end{figure}

PyLith~\cite{pylith-web-page,AagaardKnepleyWilliams13} is a finite element code with a primary focus on modeling interseismic and coseismic deformation of Earth’s crust and upper mantle.
PyLith supports 2D and 3D static, quasistatic (neglecting inertia), and dynamic (including inertia) formulations of the governing equations, which can be elasticity, incompressible elasticity, or poroelasticity, with a variety of elastic and viscoelastic bulk rheologies.
Dirichlet and Neumann boundary conditions are supported, as well as absorbing boundaries.
PyLith is able to model faults as dislocations along interior interfaces, both prescribed slip (kinematic rupture) and spontaneous rupture (fault friction).
The code is written in C++ and Python. % and uses MPI for parallel processing.
It leverages PETSc~\cite{petsc-user-ref,petsc-web-page} for finite element data structures and operations as well as linear and nonlinear solvers.
This subduction problem involves a quasi-static simulation for interseismic deformation, meaning that the inertial effects are neglected.
We prescribe aseismic slip (creep) on the bottom of the slab and the deeper portion of the top of the slab; the shallow portion of the top of the slab remains locked.
We use linear Maxwell viscoelastic bulk rheologies in the mantle and deeper part of the slab.
Fig.~\ref{fig:subduction-3d-step03-diagram} shows the boundary conditions on the domain.

The equation of the conservation of momentum on the fault interface reduces to
\begin{equation}
  \int_{\Gamma_{\text{f}^+}} \Sigmab \cdot \mathbf{n} + \mathbf{\lambda} \, d\Gamma + \int_{\Gamma_{\text{f}^-}} \Sigmab \cdot \mathbf{n} - \mathbf{\lambda} \, d\Gamma = 0,
\end{equation}
in which $\Sigmab$ is the stress tensor, $\Gamma_\text{f}$ is the fault surface, and $\mathbf{\lambda}$ is the Lagrange multiplier enforcing the prescribed slip along the fault. We enforce this equation on each portion of the fault interface along with our prescribed slip constraint, which leads to
\begin{align}
  \Sigmab \cdot \mathbf{n} + \mathbf{\lambda} &= \mathbf{0} \text{ on } \Gamma_{\text{f}^+}, \\
  \Sigmab \cdot \mathbf{n} - \mathbf{\lambda} &= \mathbf{0} \text{ on } \Gamma_{\text{f}^-}, \\
  \ub^+ - \ub^- - \mathbf{d}(\xb,t) &= \mathbf{0},
\end{align}
in which $\ub$ is the displacement vector and $\mathbf{d}$ is the slip.
Our solution vector consists of both displacements and Lagrange multipliers, and the strong form for the system of equations is
\begin{align}
  % Solution
  \mathbf{s}^T &= \left( \ub \quad \mathbf{\lambda} \right)^T \\
  % Elasticity
  \mathbf{f}(\xb,t) + \boldsymbol{\nabla} \cdot \Sigmab(\ub) &= \mathbf{0} \text{ in }\Omega, \\
  % Neumann
  \Sigmab \cdot \mathbf{n} &= \mathbf{\tau}(\xb,t) \text{ on }\Gamma_{\mathbf{\tau}}, \\
  % Dirichlet
  \ub &= \ub_0(\xb,t) \text{ on }\Gamma_u, \\
  % Prescribed slip
  \ub^+ - \ub^- - \mathbf{d}(\xb,t) &= \mathbf{0} \text{ on }\Gamma_{\text{f}},  \\
  \Sigmab \cdot \mathbf{n} &= -\mathbf{\lambda}(\xb,t) \text{ on }\Gamma_{\text{f}^+}, \\
  \Sigmab \cdot \mathbf{n} &= +\mathbf{\lambda}(\xb,t) \text{ on }\Gamma_{\text{f}^-},
\end{align}
in which $\mathbf{\tau}$ is the traction and $\mathbf{f}$ is the body force.
PyLith the problem using a linear finite element for the displacement, and solve using a fully implicit timestepper (backwards Euler).
A complete description of this problem can be found in the PyLith User's Manual~\cite{PyLithManual}, or in the \texttt{examples/subduction-3d} directory of the source tree.

%% TODO, at 'we used the lower tridiagonal'... ``Need to say that this is with PCFIELDSPLIT?''
The linear system is solved using a Schur complement formulation.
The displacement and Lagrange multiplier blocks were solved directly using sparse direct factorization, in order to remove the complication of the interaction of the subsolver tolerance with the outer loop convergence.
However, these systems are often solved with algebraic multigrid.
We used the lower tridiagonal variant of the Schur factorization (the \texttt{lower} variant in PETSc), and approximated the elastic block with its diagonal in order to form a preconditioning matrix for the Lagrange block (the \texttt{selfp} variant in PETSc).
The system used both a relative and absolute tolerance of $10^{-12}$.

As the fault system is not positive definite, we do not report results for the second Fischer initial guess method (which requires a symmetric positive definite matrix) for this benchmark.
For a small number of vectors, both POD and \WING outperform the original Fischer algorithm in terms of iterations.
As the number of vectors is increased, all variants tend toward the same number of iterates.
However, this is still almost a factor of $2.5$ fewer iterates than the initial solver currently used in PyLith.
Using runtime as a metric, POD is the clear winner.

\begin{figure}
  \centering
  \begin{tabular}{c c}
  \includegraphics[width=0.45\linewidth]{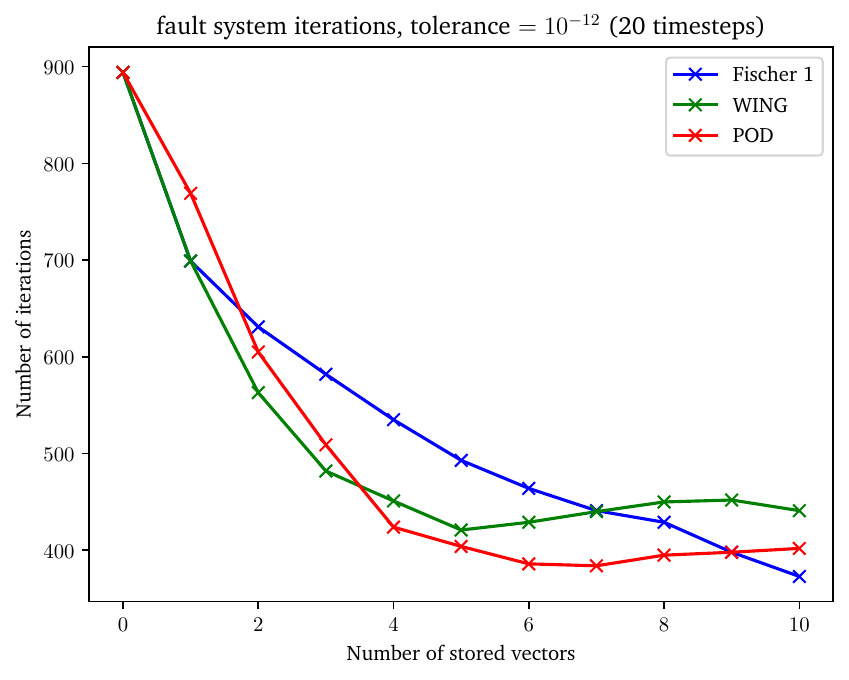} &
  \includegraphics[width=0.45\linewidth]{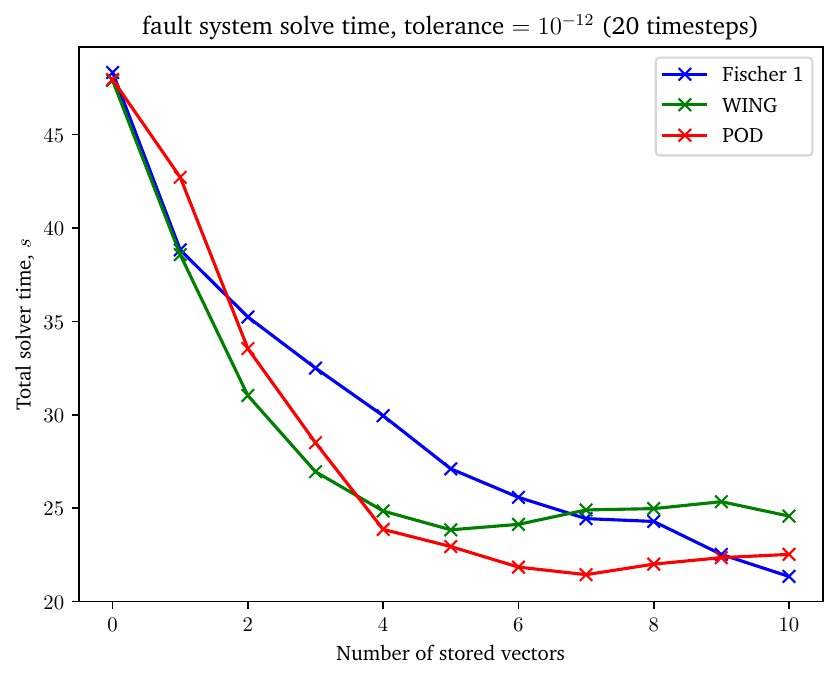} \\
  \end{tabular}

  \caption{Solver iteration numbers and times for the 3D subduction experiment for finite element system from PyLith.}

  \label{fig:pylith}
\end{figure}

\subsection{Mantle Convection}
\label{subsec:mantle-convection}
This example is based on the step-32 tutorial program which is included in the deal.II finite element library~\cite{dealII96}.
deal.II is a parallel finite element library written in C++ with support for distributed linear algebra using PETSc and the \emph{HYPRE}~\cite{hypre-web-page} preconditioner library.
A thorough description of this program's governing equations, modeling assumptions, and numerical methods is available in Kronbichler et al.~\cite{step-32}.

\begin{figure}
\centering
% scale is 973 to 4273
\begin{tabular}{c c}
\includegraphics[width=0.4\textwidth]{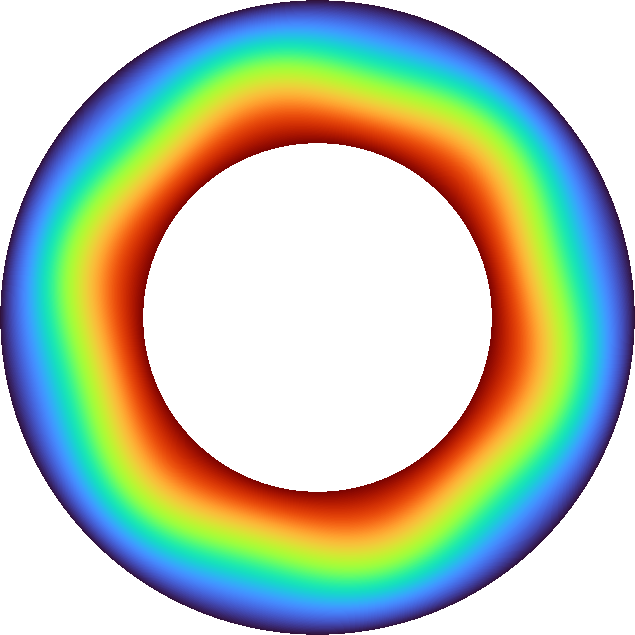} \hspace{0.05\textwidth} &
\includegraphics[width=0.4\textwidth]{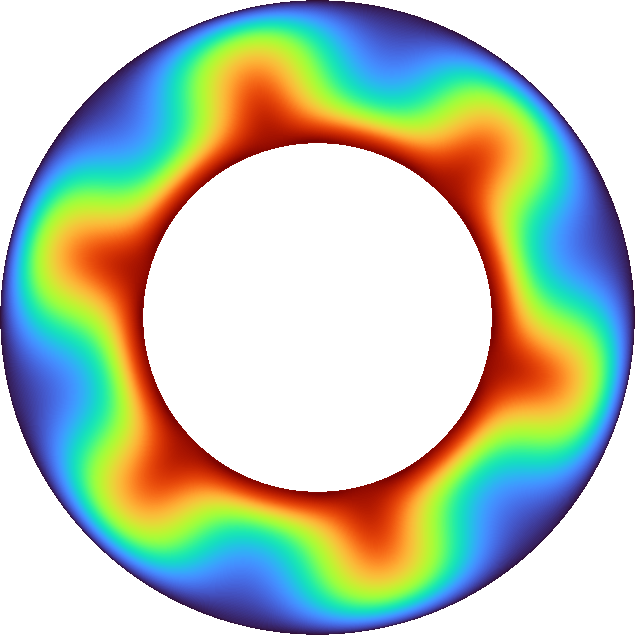} \hspace{0.05\textwidth}
\end{tabular}

% This spacing is eye-balled
\vspace{0.05\textwidth}
\includegraphics[width=0.4\textwidth]{Plots/colorbar.png}\\
$973 \text{ K}$\hspace{0.32\linewidth}$4273 \text{ K}$

\caption{Temperature field for the mantle convection benchmark at the initial condition and after $200$ time steps.}
\end{figure}

This program uses the finite element method to solve the coupled mantle convection equations
\begin{align}
  -\nabla \cdot (2 \nu \varepsilon(\ub)) + \nabla p &= \rho(T) \boldsymbol{g}, \\
  \nabla \cdot \ub &= 0, \\
  \dfrac{\partial T}{\partial t} + \ub \cdot \nabla T - \nabla \cdot \kappa \nabla T &= \gamma, \\
  \rho(T) &= \rho_0 (1 - \beta (T - T_0)),
\end{align}
in which $\ub$ is the fluid velocity, $\varepsilon(\ub)$ is the symmetrized gradient of the velocity, $\nu$ is the viscosity, $p$ is the total pressure, $\rho$ is the fluid density, $T$ is the fluid temperature, $g$ is the gravity, $\kappa$ is the thermal diffusivity, and $\gamma$ is the rate of internal heating (which depends on $\ub$).
Density is modeled as a linear function of temperature with reference temperature and density $(T_0, \rho_0)$.

step-32 uses Dirichlet boundary conditions of $T = 4273 \text{ K}$ along the inner boundary and $T = 973 \text{ K}$ along the outer boundary.
For simplicity, we use zero velocity along the inner boundary and zero normal velocity along the outer boundary.
These boundary conditions are relatively simple and suffice for this benchmark.
More sophisticated models are possible and these tradeoffs are discussed in Kronbichler et al.~\cite{step-32}.
Unlike the original publication, for the sake of simplicity we do not use adaptive mesh refinement (AMR) as slight changes to solution vectors, even within a specified numerical tolerance, can result in slightly different cells being refined or coarsened in time, which in turn alters the number of degrees of freedom.
Hence, by disabling AMR, we guarantee that we are using the same linear system with all initial guess algorithms.

The resulting linear system for velocity and pressure is the standard Stokes system
\begin{equation}
  \begin{pmatrix}
    A & B^T \\
    B & 0
  \end{pmatrix}
  \begin{pmatrix}
    U \\
    P
  \end{pmatrix}
  =
  \begin{pmatrix}
    F_U \\
    0
  \end{pmatrix}
\end{equation}
in which $A$ is a discretization of the Laplacian and $B^T$ is a discretization of the gradient.
We solve this system with the FGMRES algorithm and a Schur complement preconditioner with a relative tolerance stopping criterion of $10^{-8}$.
This system's preconditioner contains two conjugate gradient solves.
The first approximately solves for the pressure by replacing the Schur complement matrix $B A^{-1} B^T$ with a pressure mass matrix (and corresponding Jacobi preconditioner).
The second solve uses an algebraic multigrid (AMG) preconditioner, implemented by \emph{HYPRE}~\cite{hypre-web-page}, to provide an approximate solver for $A U^* = B^T P^* + U$, in which $P^*$ is the pressure calculated by solving the pressure mass matrix and $U$ is the provided velocity iterate.
Using the initial guess routines on the outermost FGMRES solver sufficiently demonstrates the relative efficacy of the three considered methods, so we do not use them with the preconditioner's two linear solvers.

The time derivative for the temperature convection equation is discretized (after the first timestep) with BDF-2.
The resulting discretization handles the diffusion term implicitly and the convection term explicitly.
Hence the equation for evolving the temperature is
\begin{align}
  \nonumber
  \dfrac{3}{2} T^n - k \nabla \cdot \kappa \nabla T^n
  &= 2 T^{n - 1} - \dfrac{1}{2} T^{n - 2}                                     \\
  &+ k \nabla \cdot \left[\nu_\alpha \nabla (2 T^{n - 1} - T^{n - 2}) \right] \\
  \nonumber
  &- k (2 \ub^{n - 1} - \ub^{n - 2}) \cdot \nabla (2 T^{n - 1} - T^{n - 2})
  + k\gamma
\end{align}
in which $n > 1$ is the timestep number, $k$ is the timestep, and $\nu_\alpha$ is a stabilization parameter discussed at length in Kronbichler et al.~\cite{step-32}.
The resulting linear system, therefore, is a mass matrix plus the time step times a Laplace matrix.
This is, in practice, sufficiently well-conditioned that we solve it with the conjugate gradient algorithm and a Jacobi preconditioner to a relative tolerance of $10^{-12}$.
For this particular problem the linear solve for the temperature is an order of magnitude less expensive than the Stokes solve, so the overall method does not benefit much from more advanced preconditioners.

\begin{figure}
  \centering
  \begin{tabular}{c c}
  \includegraphics[width=0.45\linewidth]{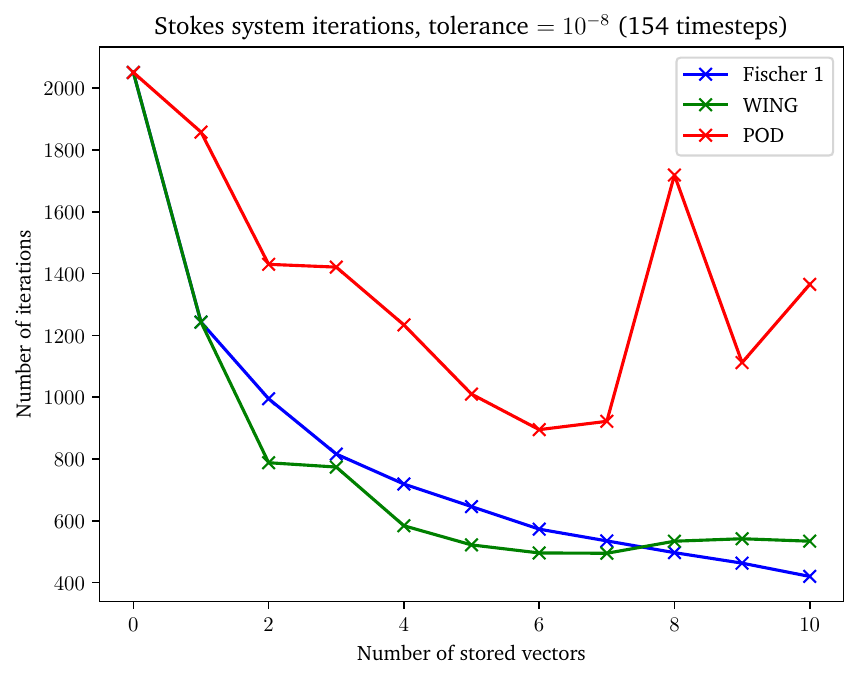} &
  \includegraphics[width=0.45\linewidth]{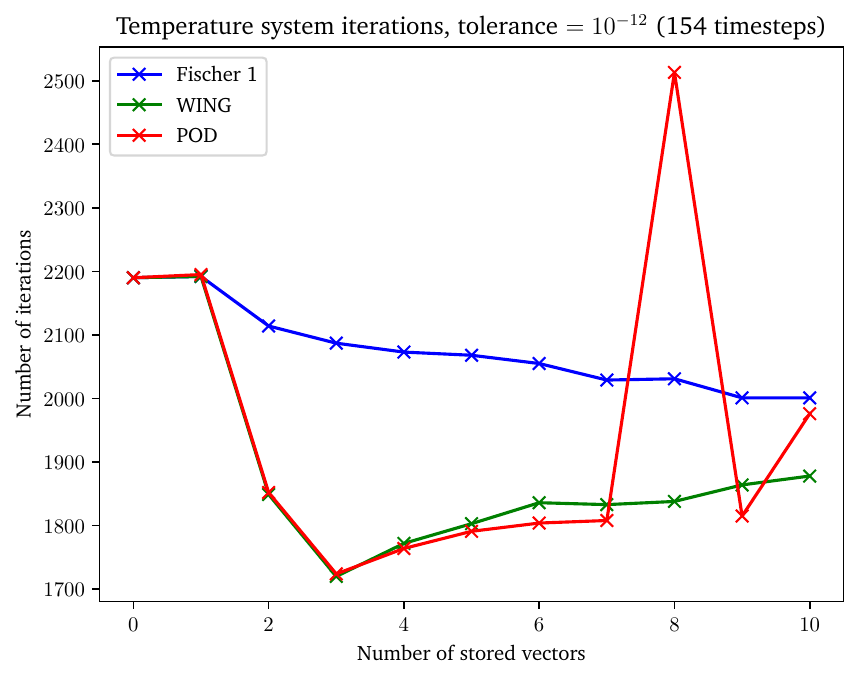} \\
  \end{tabular}

  \caption{Iteration counts for the mantle convection benchmark.}

  \label{fig:step-32-iterations}
\end{figure}

\begin{figure}
  \centering
  \begin{tabular}{c c}
  \includegraphics[width=0.45\linewidth]{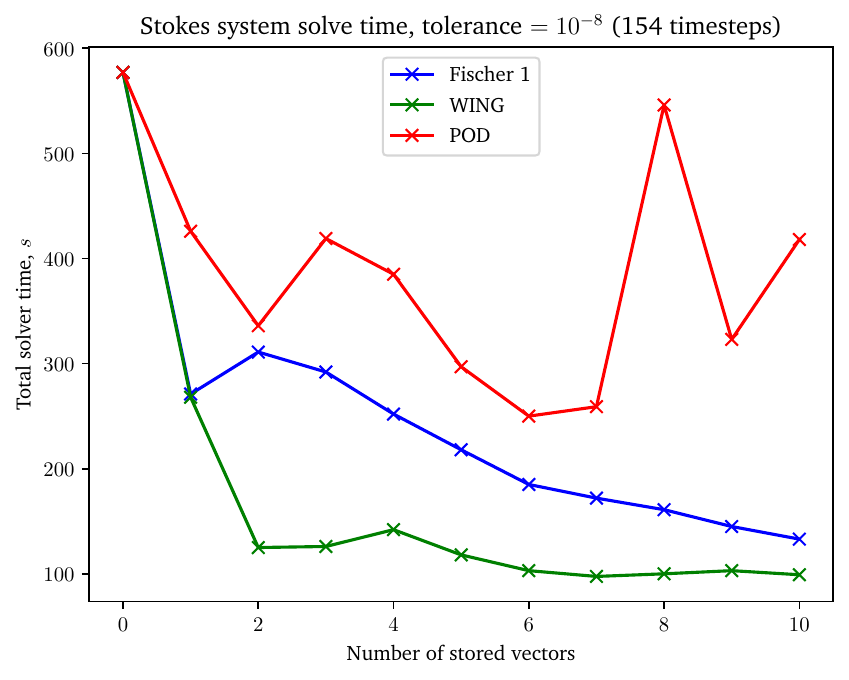} &
  \includegraphics[width=0.45\linewidth]{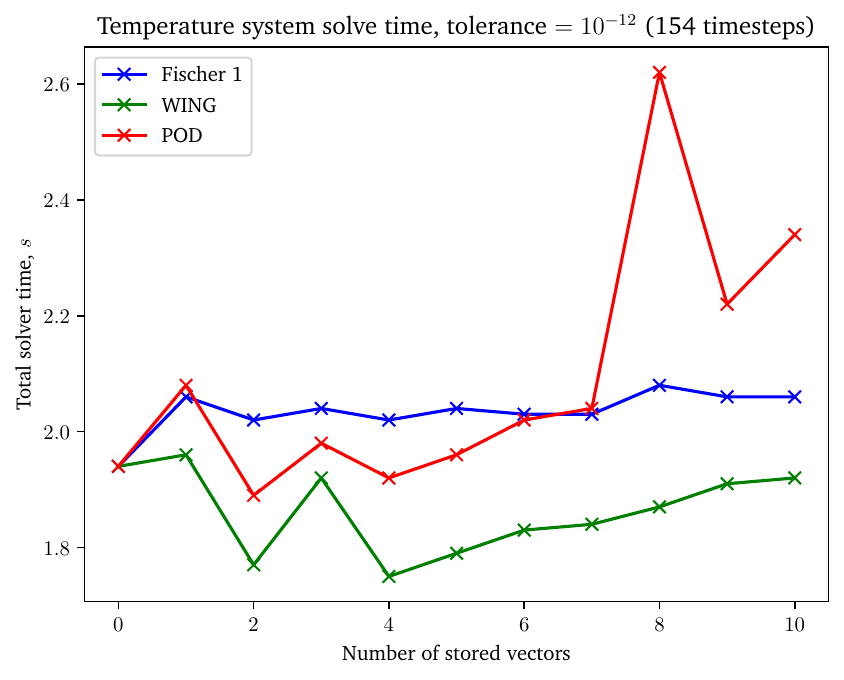} \\
  \end{tabular}

  \caption{Solver times for the mantle convection benchmark.}

  \label{fig:step-32-time}
\end{figure}

We do not report results for the second Fischer initial guess method for this benchmark because that method requires a symmetric positive definite linear operator, a property the Stokes system does not have.
All three initial guess techniques succeed in lowering the total number of iterations and total solver time.
However, the temperature system is sufficiently well conditioned (and using the previous solution vector is a good enough initial guess) that it does not benefit much from any of the initial guess algorithms: at best they lower the total solver time by about 10\%.
In contrast, the Stokes system benefits significantly from improved initial guesses, which lower the total solver time and number of iterations by about 80\%.
The POD initial guess routine performs significantly worse, in both solver time and number of iterations, with eight vectors than any other number of vectors.
This breakdown does not appear to have a single direct cause.
This initial guess algorithm does not have any guarantees on its performance so, to the best of our knowledge, this seems to simply be a case in which the generated POD vectors do not accurately represent the solution of the system.

\section{Discussion and Conclusion}
This work introduces a new initial guess method based on projection without orthogonalization.
Unlike previous projection methods, such as those devised by Fischer~\cite{fischer1998projection}, this method uses sets of possibly linearly dependent vectors rather than bases of linear spaces.
This new method, the \WING algorithm, is summarized in Algorithms~\ref{alg:form} and \ref{alg:update}.
\WING is essentially optimal in the sense that it forms initial guesses by projecting a right-hand side onto a specified set of vectors with a minimal number of dot products and scaled vector additions while using no matrix-vector products.

We examine the efficacy of the \WING algorithm for both lowering the number of linear solver iterations and reducing the overall time spent in linear solvers.
To this end, we selected a diverse set of benchmarks corresponding to fluid-structure interaction, crustal deformation modeling, and mantle convection.
Furthermore, we ran the same set of benchmarks for three other initial guess algorithms: Fischer's two methods~\cite{fischer1998projection} and the POD-based method of Zampini~\cite{zampini2010non}.
To the best of our knowledge, these numerical experiments are the first direct comparison of all three methods for the same benchmarks in a single work.
For very strict numerical tolerances, Fischer's methods tended to have the fastest time to solution. In all other cases, \WING was either the fastest method or competitive with the fastest method.
In particular, \WING tended to work best with only two or three stored solution-right-hand side vector pairs.

% main text bibliography:
\renewcommand\refname{REFERENCES}
\bibliographystyle{plain}
\bibliography{bib}

\end{document}